\documentclass[a4paper,10pt,fleqn]{article}
\usepackage{amsbsy}
\usepackage{amscd}
\usepackage{amsgen}
\usepackage{amsmath}
\usepackage{amssymb}
\usepackage{amstext}
\usepackage{amsthm}
\usepackage{upref}
\usepackage{graphicx}

\newtheorem{theoremnonum}{Theorem}

\newtheorem{theorem}{Theorem}[section]
\newtheorem{lemma}[theorem]{Lemma}

\newtheorem{proposition}[theorem]{Proposition}

\theoremstyle{definition}

\newtheorem{definition}{Definition}[section]

\newtheorem{remark}{Remark}[section]

\title{Scharlemann-Thompson untelescoping of Heegaard splittings is 
finer than Casson-Gordon's}
\author{Tsuyoshi Kobayashi}
\date{}
\begin{document}
\maketitle

\section{Introduction}
Let $H_1 \cup_P H_2$ be a Heegaard splitting of a closed 3-manifold 
$M$, 
i.e., 
$H_i$ $(i=1,2)$ is a handlebody in $M$ such that 
$H_1 \cup H_2 = M$, $H_1 \cap H_2 = \partial H_1 = \partial H_2 = P$. 
In \cite{ST}, M.Scharlemann, and A.Thompson had 
introduced a process for spreading 
$H_1 \cup_P H_2$ into a \lq\lq thinner\rq\rq\/ presentation. 
The idea was polished to show that if the original Heegaard splitting 
is irreducible, then we can spread it into a series 
$(A_1 \cup_{P_1} B_1) \cup \cdots \cup (A_n \cup_{P_n} B_n)$ 
such that each $A_i \cup_{P_i} B_i$ is a strongly irreducible 
Heegaard splitting. 
In this paper, we call this series of strongly irreducible Heegaard splittings a 
{\it Scharlemann-Thompson untelescoping} (or {\it S-T untelescoping}) 
of $H_1 \cup_P H_2$. 
On the other hand, preceding \cite{ST}, A.Casson, and C.Gordon \cite{CG} 
had proved that if $H_1 \cup_P H_2$ is weakly reducible and not reducible, 
then there exists an incompressible surface of positive genus in $M$. 
This result is proved by using the following argument. 

Let $\Delta = \Delta_1 \cup \Delta_2$ be a weakly 
reducing collection of compressing disks for $P$ 
(for the definitions of the terms, see section~2). 
Then $P( \Delta )$ denotes the surface obtained from $P$ by compressing 
along $\Delta$. 
Let $\hat{P}( \Delta )$ be the surface obtained from $P( \Delta )$ 
by discarding the components that are contained in either $H_1$ or $H_2$. 
%of $\hat{P}( \Delta )$ has positive genus. 
Suppose that $\Delta$ has minimal complexity 
(for the definition of the complexity, see section~4). 
Then we can show that the irreducibility of $H_1 \cup_P H_2$ implies 
that no component of $\hat{P}( \Delta )$ is a 2-sphere. 
Then, by using a relative version of Haken's theorem \cite{Ha}, 
we can show that $\hat{P}( \Delta )$ is incompressible. 

With adopting the above notations, we will see, in section~4, 
that the closure of each component of $M-\hat{P}( \Delta )$ naturally 
inherits a Heegaard splitting from $H_1 \cup_P H_2$ 
if $\hat{P}( \Delta )$ contains no 2-sphere component. 
Hence we obtain a series of Heegaard splitings, say 
$(C_1 \cup_{Q_1} D_1) \cup \cdots \cup (C_m \cup_{Q_m} D_m)$. 
If $\hat{P}( \Delta )$ is incompressible, then this series is called a 
{\it Casson-Gordon untelescoping} (or {\it C-G untelescoping}) 
of $H_1 \cup_P H_2$. 
Then we will also see that C-G untelescoping of $H_1 \cup_P H_2$ 
can be regarded as one that 
appears in a process for obtaining S-T untelescoping from $H_1 \cup_P H_2$
(Remark~\ref{C-G induces S-T}). 

We note that these two untelescopings are used in many articles and 
equally usefull. 
For example, C-G untelescoping was used 
by M.Boileau, and J.-P.Otal \cite{BO} for studying 
Heegaard splittings of the 3-dimensional torus, 
by J.Schultens \cite{Schu'} for studying Heegaard splittings 
of $\text{(surface)} \times S^1$, 
by M.Lustig, and Y.Moriah \cite{LM} for studying the exteriors of 
wide knots and links, 
and by the author  \cite{Ko2} for studying the Heegaard splittings 
of the exteriors of two bridge knots. 
S-T untelescoping was used, for example, 
by Scharlemann-Schultens \cite{SS}, Schultens\cite{Schu}, Morimoto\cite{Mo}, and 
Morimoto-Schultens \cite{MS}for studying the Heegaard splittings 
of the exteriors of non-prime knots. 
However, it seems that it is not known that whether these two concepts 
are the same one or not. 

Hence, it is natural to ask:

\begin{quote}
Question. Are these two untelescopings essentially the same ?
\end{quote}

Since C-G untelescoping of $H_1 \cup_P H_2$ can be regarded 
as an untelescoping that appeares in a process for obtaining S-T untelescoping 
from $H_1 \cup_P H_2$, the above question can be strengthened 
as in the following form. 

\begin{quote}
Question$'$. Is S-T untelescoping essentially finer than C-G untelescoping ?
\end{quote}

The purpose of this paper is to show that the answer to Question$'$ is 
positive.

\begin{theoremnonum}\label{main theorem}
There exist infinitely many closed, orientable, Heegaard genus 4 3-manifolds 
such that each 3-manifold $M$ admits a genus 4 Heegaard splitting 
$V \cup_P W$ with the following properties. 

\begin{enumerate}

\item 
There is a S-T untelescoping of $V \cup_P W$ which decomposes $M$ into 
three pieces, say 
$(V_1 \cup_{P_1} W_1) \cup (V_2 \cup_{P_2} W_2) \cup (V_3 \cup_{P_3} W_3)$. 

\item 
The Heegaard splitting $V \cup_P W$ is decomposed into exactly two pieces 
by any C-G untelescoping. 

\item
There is a C-G untelescoping 
$(V_1 \cup_{P_1} W_1) \cup (V_2' \cup_{P_2'} W_2')$ of 
$V \cup_P W$, such that $V_1 \cup_{P_1} W_1$ is the Heegaard 
splitting that appeared in the above 1, and that 
$V_2' \cup_{P_2'} W_2'$ is weakly reducible. 
Moreover, $(V_2 \cup_{P_2} W_2) \cup (V_3 \cup_{P_3} W_3)$ is 
a C-G untelescoping of $V_2' \cup_{P_2'} W_2'$, 
where $V_2 \cup_{P_2} W_2$ and  $V_3 \cup_{P_3} W_3$ 
are Heegaard splittings that appeared in the above 1. 
\end{enumerate}
\end{theoremnonum}

\section{Preliminaries}

Throughout this paper, 
we work in the piecewise linear category. 
For a submanifold $H$ of a manifold $M$, 
$N(H,M)$ denotes a regular neighborhood of $H$ in $M$. 
When $M$ is well understood, we often abbreviate $N(H,M)$ to $N(H)$. 
Let $N$ be a manifold embedded in a manifold $M$
with dim$N=$dim$M$.
Then $\text{Fr}_M N$ denotes the frontier of $N$ in $M$.
For the definitions of standard terms in 
3-dimensional topology, 
we refer to \cite{He} or \cite{Ja}. 

A 3-manifold $C$ is a {\it compression body} if there exists a compact 
connected closed surface $F$ 
such that $C$ is obtained 
from $F \times [0,1]$ by attaching 2-handles along mutually disjoint 
simple closed curves in $F \times \{ 1 \}$ and capping off the resulting 
2-sphere boundary components which are disjoint from $F \times \{ 0 \}$ 
by 3-handles. 
The subsurface of $\partial C$ corresponding to $F \times \{ 0 \}$ 
is denoted by $\partial_+C$. 
Then $\partial_-C$ denotes the subsurface 
$\partial C - \partial_+C$ 
of $\partial C$. 
A compression body $C$ is said to be {\it trivial} if 
$C$ is homeomorphic to $F \times [0,1]$ with $\partial_-C$ corresponding 
to $F \times \{ 0 \}$. 
A compression body $C$ is called a {\it handlebody} if 
$\partial_-C = \emptyset$. 
A compressing disk $D (\subset C)$ of $\partial_+ C$ is called 
a {\it meridian disk} of the compression body $C$. 

\begin{remark}\label{remark of compression body}
The following properties are known for compression bodies. 

\begin{enumerate}

\item 
The compression bodies are irreducible. 

\item 
By extending the cores of the 2-handles in the definition of 
the compression body $C$ vertically to $F \times [0,1]$, 
we obtain a union of mutually disjoint meridian disks 
${\cal D}$ of $C$ such that the manifold 
obtained from $C$ by cutting along ${\cal D}$ is homeomorphic 
to a union of $\partial_- C \times [0,1]$ and some 
(possibly empty) 3-balls. 
This gives a {\it dual description} of compression bodies. 
That is, a connected 3-manifold $C$ is a compression body 
if there exists a compact (not necessarily connected) closed surface 
${\cal F}$ without 2-sphere components and a union of 
(possibly empty) 3-balls ${\cal B}$ such that $C$ is obtained from 
${\cal F} \times [0,1] \cup {\cal B}$ by attaching 1-handles to 
${\cal F} \times \{ 0 \} \cup \partial {\cal B}$. 
We note that $\partial_-C$ is the surface corresponding to 
${\cal F} \times \{ 1 \}$.

\end{enumerate}
\end{remark}

Let $N$ be a cobordism between two closed surfaces 
$F_1$, $F_2$ (possibly $F_1 = \emptyset$ or $F_2 = \emptyset$), 
i.e., 
$F_1 \cup F_2$ is a partition of the components of $\partial N$. 

\begin{definition}\label{Heegaard splitting}
We say that $C_1 \cup_P C_2$ (or $C_1 \cup C_2)$ 
is a {\it Heegaard splitting} 
of $(N, F_1, F_2)$ (or simply, $N$) if it satisfies the following 
conditions. 

\begin{enumerate}
\item 
$C_i$ $(i=1,2)$ is a compression body in $N$ such that 
$\partial_- C_i = F_i$, 

\item 
$C_1 \cup C_2 = N$, and 

\item 
$C_1 \cap C_2 = \partial_+ C_1 = \partial_+ C_2 = P$.

\end{enumerate}

The surface $P$ is called a {\it Heegaard surface} of 
$(N, F_1, F_2)$ (or, $N$). 
The genus of 
$P$ is called the {\it genus} of the Heegaard splitting. 
\end{definition}

\begin{definition}\label{reducible Heegaard splitting}

\ 
\begin{enumerate}

\item 
A Heegaard splitting $C_1 \cup_P C_2$ is {\it reducible} 
if there exist meridian disks $D_1$, $D_2$ of 
the compression bodies $C_1$, $C_2$ 
respectively such that 
$\partial D_1= \partial D_2$ 

\item 
A Heegaard splitting $C_1 \cup_P C_2$ is {\it weakly reducible} 
if there exist meridian disks $D_1$, $D_2$ of 
the compression bodies $C_1$, $C_2$ 
respectively such that 
$\partial D_1 \cap  \partial D_2 = \emptyset$. 
If $C_1 \cup_P C_2$ is not weakly reducible, then it is called 
{\it strongly irreducible}.

\item 
A Heegaard splitting $C_1 \cup_P C_2$ is {\it trivial} 
if either $C_1$ or $C_2$ is a trivial compression body. 
\end{enumerate}
\end{definition}

\section{Scharlemann-Thompson untelescoping}

Let $C_1 \cup_P C_2$ be a Heegaard splitting of $(M, F_1, F_2)$. 
By 2 of Remark~\ref{remark of compression body}, 
we see that $C_1$ is obtained from 
$F_1 \times [0,1] \cup 0$-handles by attaching 1-handles. 
Recall that $C_2$ is obtained from 
$\partial_+ C_2 \times [0,1]$
by attaching 2-handles, and 3-handles. 
Then, by using an isotopy which pushes 
$\partial_+ C_2 \times [0,1]$ out of $C_2$, we identify 
$\partial_+ C_2 \times [0,1]$ with $N( \partial_+ C_1, C_1 )$. 
This identification together with the above handles gives the 
following handle decomposition of $M$. 

$$M=F_1 \times [0,1] \cup ( 0\text{-handles}) \cup ( 1\text{-handles}) \cup 
( 2\text{-handles}) \cup ( 3\text{-handles})$$

We note that there are huge variety of ways for giving 
handle decompostions for $H_1$, $H_2$. 
Suppose that: 

\begin{quote}
there exists a proper subset of the 
$0\text{-handles} \cup 1\text{-handles}$ such that 
some subset the 
$2\text{-handles} \cup 3\text{-handles}$ 
do not intersect the 0-handles and 1-handles at all. 
\end{quote}

\noindent
Then we can arrange the order of the handles non-trivially to obtain 
submanifolds 
$N_1, \dots , N_n$ 
such that 

\noindent
$N_1 = 
(F_1^{(1)} \times [0,1]) \cup 
( 0\text{-handles}) \cup ( 1\text{-handles}) \cup 
( 2\text{-handles}) \cup ( 3\text{-handles}), $

\noindent
$N_2 = N_1 \cup 
(F_1^{(2)} \times [0,1]) \cup 
( 0\text{-handles}) \cup ( 1\text{-handles}) \cup 
( 2\text{-handles}) \cup ( 3\text{-handles}), $

$\dots$

\noindent
$N_n  = 
N_{n-1} \cup 
(F_1^{(n)} \times [0,1]) \cup 
( 0\text{-handles}) \cup ( 1\text{-handles}) \cup 
( 2\text{-handles}) \cup ( 3\text{-handles}), $

\noindent
where $F_1^{(1)} \cup F_1^{(2)} \cup \cdots \cup F_1^{(n)}$ is a 
partition of the components of $F_1$, and 
each handle is one from the handle decompositon of $M$. 
Suppose that the following properties are satisfied. 

\begin{enumerate}

\item 
$N_1$ is connected, and $N_n = M$. 

\item
At each stage $k$ $(2 \le k \le n)$, 
let $\partial_k^-$ denotes the union of the components of 
$\partial N_{k-1}$ to which the $1$-handles are attached. 
Then 
$\partial_k^- \cup (F_1^{(k)} \times [0,1]) \cup 
( 0\text{-handles}) \cup ( 1\text{-handles}) \cup 
( 2\text{-handles}) \cup ( 3\text{-handles})$ 
is connected, where these handles are those that 
appeared in the stage $k$. 

\item
Each component of 
$\partial N_1, \partial N_2, \dots , \partial N_{n-1}$ 
is not a 2-sphere. 
\end{enumerate}

Then at each stage $k$ 
let $I_k = \partial_k^- \times [0,1]$, 
and 
$C_1^{(k)} = 
I_k \cup (F_1^{(k)} \times [0,1]) 
\cup ( 0\text{-handles}) \cup ( 1\text{-handles})$, 
where 
$\partial_k^- \times \{ 0 \}$ and $\partial_k^-$ are identified, 
and (0-handles), (1-handles) are those that 
appeared in the stage $k$. 
Since each component of $\partial N_k$ is not a 2-sphere 
(above condition~3), 
we see that each component of $C_1^{(k)}$ is a compression body 
by 2 of Remark~\ref{remark of compression body}. 
Then we see that $C_1^{(k)}$ is connected (above condition~2), 
hence $C_1^{(k)}$ is a compression body. 
Hence $\partial_+ C_1^{(k)}$ is a connected surface, 
and let 
$J_k = \partial_+ C_1^{(k)} \times [0,1]$, and 
$C_2^{(k)} = J_k \cup ( 2\text{-handles}) \cup ( 3\text{-handles})$, 
where 
$\partial_+ C_1^{(k)} \times \{ 1 \}$ and $\partial_+ C_1^{(k)}$ 
are identified, and 
(2-handles), (3-handles) are those that 
appeared in the stage $k$. 
Then we see that $C_2^{(k)}$ is a compression body by the above condition~3. 
It is clear from the construction that 
we obtained a submanifold, say $R_k$, 
of $M$ with the Heegaard splitting $C_1^{(k)} \cup C_2^{(k)}$. 
Moreover it is clear that $M$ can be regarded as obtained from 
$R_1, \dots , R_n$ by identifying their boundaries. 
We call the decomposition of 
$M$ into the series of Heegaard splittings 
$( C_1^{(1)} \cup C_2^{(1)}) \cup \cdots \cup ( C_1^{(n)} \cup C_2^{(n)})$ 
an {\it untelescoping} of the Heegaard splitting 
$C_1 \cup_P C_2$. 

\begin{definition}\label{S-T untelescoping}
The above untelescoping is called a 
{\it Scharlemann-Thompson untelescoping} 
(or {\it S-T untelescoping}) 
if each Heegaard splitting 
$C_1^{(k)} \cup C_2^{(k)}$ is non-trivial, and strongly irreducible. 
\end{definition}

\begin{remark}\label{Remark of S-T untelescoping}
It is known that every irreducible Heegaard splitting of 
3-manifolds with incompressible boundary admits a S-T untelescoping 
(see, for example, \cite{S}). 
\end{remark}

\section{Casson-Gordon untelescoping}

In \cite{CG}, 
A.Casson and C.McA.Gordon proved that if a Heegaard splitting of a 
closed 3-manifold is weakly reducible, and not reducible, then 
$M$ contains an incompressible surface of positive genus. 
In this section, 
we introduce the arguments in their proof, and 
give the definition of Casson-Gordon untelescoping. 

Let $M$ be a closed, orientable 3-manifold, and 
$H_1 \cup_P H_2$ a Heegaard splitting of $M$. 
Let 
$\Gamma = \Gamma_1 \cup \Gamma_2$ 
be a weakly reducing collection of disks for $P$, 
i.e., 
$\Gamma_i$ $(i=1,2)$ is a union of mutually disjoint, 
non-empty meridian disks of $H_i$ such that 
$\Gamma_1 \cap \Gamma_2 = \emptyset$. 
Then $P( \Gamma )$ denotes the surface obtained from $P$ 
by compressing $P$ along $\Gamma$. 
Let 
$\hat{P}( \Gamma ) = P( \Gamma ) - ($the components of 
$P( \Gamma )$ 
which are contained in $H_1$ or $H_2)$. 
In \cite{CG}, Casson-Gordon proved: 

\begin{proposition}\label{CG Theorem 3.1}
Let $M$ be a closed, orientable 3-manifold, 
and $H_1 \cup_P H_2$ a Heegaard splitting of $M$. 
Suppose that $H_1 \cup_P H_2$ is weakly reducible. 
Then either 

\begin{enumerate}
\item 
$H_1 \cup_P H_2$ is reducible, or 

\item
there exists a weakly reducing collection of disks 
$\Delta$ for $P$ such that each component of 
$\hat{P}( \Delta )$ is an incompressible surface in $M$, 
which is not a 2-sphere. 
\end{enumerate}
\end{proposition}

This result is proved by using the following argument. 

\begin{quote}
In general, for a closed surface $F$, we define a complexity $c(F)$ 
of $F$ as follows. 

$$c(F) = \Sigma (1- \chi (F^i )),$$ 

\noindent 
where the sum is taken over each component $F^i$ 
of $F$ that is not a 2-sphere. 

Let $\Delta = \Delta_1 \cup \Delta_2$ be a weakly reducing 
collection of disks for $P$ such that $c( \hat{P}( \Delta ) )$ is minimal 
among the weakly reducing collection of disks of $P$. 
We can show that if $\hat{P}( \Delta )$ contains a 2-sphere component, 
then $H_1 \cup_P H_2$ is reducible, and this gives conclusion~1. 
(We note that, for the proof of this assertion, the maximality of 
$\Delta$ is not necessary. )
If no component of $\hat{P}( \Delta )$ is a 2-sphere, 
then the minimality of $c( \hat{P}( \Delta ) )$ 
together with a relative version of 
Haken's lemma shows that each component of $\hat{P}( \Delta )$ is 
incompressible, and this gives conclusion~2, and 
this completes the proof of the proposition. 
\end{quote}

Now, we introduce some terminologies. 
Let $H_1 \cup_P H_2$ be a Heegaard splitting of a closed $M$, 
$\Delta = \Delta_1 \cup \Delta_2$ a weakly reducing collection of 
disks for $P$. 
Suppose that $H_1 \cup_P H_2$ is not reducible, 
hence no component of $\hat{P}( \Delta )$ is a 2-sphere. 
Let 
$M_1, \dots , M_n$ be the closures of the components of 
$M - \hat{P}( \Delta )$. 
Let 
$M_{j,i} = M_j \cap H_i$ 
$(j=1, \dots , n, i = 1,2)$. 

\begin{lemma}\label{small and big decomposition}
For each $j$, we have either one of the following. 

\begin{enumerate}

\item 
$M_{j,2} \cap P \subset \text{Int} (M_{j,1} \cap P)$, 
and 
$M_{j,1}$ is connected. 

\item 
$M_{j,1} \cap P \subset \text{Int} (M_{j,2} \cap P)$, 
and 
$M_{j,2}$ is connected. 

\end{enumerate}
\end{lemma}

\begin{proof}
Recall that $P( \Delta )$ is a surface obtained from $P$ by compressing 
along $\Delta$. 
Let 
$M_1', \dots, M_m'$ be the closures of the components of 
$M - P( \Delta )$. 
Then each $M_{\ell}'$ is either one of the following types. 

\noindent
(1) $M_{\ell}'$ is a component of 
$\text{cl}(H_i - N( \Delta_i, H_i ))$ 
$(i=1 \text{ or } 2)$. 

\noindent
(2) There is a component $E$ of 
$\text{cl}(H_i - N( \Delta_i, H_i ))$ 
$(i=1 \text{ or } 2)$ 
such that $M_{\ell}'$ is the union of $E$ and the components of 
$N( \Delta_{3-i}, H_{3-i} )$ intersecting $E$. 

By the definition of $\hat{P}( \Delta )$, 
we have: 

\medskip
\noindent
Claim 1. 
Each type~(1) 
$M_{\ell}' (\subset H_i)$ is amalgamated to the adjacent component 
in $H_{3-i}$ for obtaining some $M_k$. 

\medskip

Let $E_i$ $(i=1$ or $2)$ be a component of of 
$\text{cl}(H_i - N( \Delta_i, H_i ))$. 
By the definition, we immediately have: 

\medskip
\noindent
Claim 2. 
If $E_i \cap \Delta_{3-i} = \emptyset$, 
then $E_i$ is a type~(1) component of $M-P( \Delta )$. 

\medskip
Suppose 
$E_i \cap N(\Delta_{3-i}, H_i) \ne \emptyset$. 
Let 
$\tilde{\Delta}_{3-i}$ be the union  of the components of $\Delta_{3-i}$ 
intersecting $E_i$. 
Let 
$\tilde{E}_{3-i}$ be the union  of the components $G$ of 
$\text{cl}(H_{3-i} - N( \Delta_{3-i}, H_{3-i} ))$ 
such that 
$G \cap P \subset E_i$. 
By the definition of $\tilde{\Delta}_{3-i}$ and Claim~2, 
we see that each component of 
$\tilde{E}_{3-i}$ is type~(1). 
Moreover it is clear that $\tilde{E}_{3-i}$ is the components 
that are amalgamated to $E_i$ as in Claim~1. 
Let 
$\tilde{M} = E_i \cup N( \tilde{\Delta}_{3-i}) \cup \tilde{E}_{3-i}$. 
Since  $\tilde{E}_{3-i}$ is already amalgamated within $\tilde{M}$, 
Claim~2 shows that no component of 
$\partial \tilde{M}$ is contained in $H_1$ or $H_2$. 
Hence $\tilde{M} = M_j$ for some $j$. 
We easily see, conversely, that each $M_j$ is obtained from some 
component of $M_{\ell}' \cap H_i$ $(i=1$ or $2)$ 
as in the above manner. 
Note that these give Lemma~\ref{small and big decomposition} with regarding 
$E_i$ as $M_{j,1}$ (conclusion~1), or $M_{j,2}$ (conclusion~2). 
\end{proof}

\begin{definition}\label{small and big components}
We call the components $M_{j,i}$ $(i=1, 2)$ 
satisfying the conclusion $i$ of 
Lemma~\ref{small and big decomposition} {\it big}. 
If $M_{j,i}$ is not big, then each component of 
$M_{j,i}$ is called {\it small}. 
\end{definition}

Note that the closure of each component of 
$H_i - \hat{P}( \Delta )$ $(i=1,2)$ is a component of 
some $M_{j,i}$. 
Hence it is either big or small. 

\begin{lemma}\label{small and big are alternate}
Big components, and small components of 
$H_i - \hat{P}( \Delta )$ appear alternately in $H_i$, 
i.e., 
no pair of big components are adjacent in $H_i$, 
and no pair of small components are adjacent in $H_i$. 
\end{lemma}

\begin{proof}
Suppose that $M_{j,i}$ is a big component, 
and $E_i$ the closure of a component of 
$H_i - \hat{P}( \Delta )$, which is adjacent to $M_{j,i}$. 
Let $E_{3-i}$ be the closure of a component of 
$H_{3-i} - \hat{P}( \Delta )$ 
such that 
$E_{3-i} \cap E_i \ne \emptyset$. 
Since $E_{3-i} \cap M_{j,i} \ne \emptyset$, 
and 
$E_{3-i} \cap P$ 
is not contained in $M_{j,i}$ 
(see the Proof of Lemma~\ref{small and big decomposition}), 
we see that $E_{3-i}$ is big. 
This shows that $E_i$ is small. 

Suppose that 
$G_i (\subset H_i)$ is a small component,
and $E_i$ the closure of a component of 
$H_i - \hat{P}( \Delta )$, which is adjacent to $G_i$. 
Let $M_{j, 3-i}$ be the big component intersecting $G_i$. 
Then we have: 
$(E_i \cap P) \cap M_{j, 3-i} \ne \emptyset$, and :
$E_i \cap P$ is not contained in 
$M_{j, 3-i}$ (see the Proof of Lemma~\ref{small and big decomposition}). 
These show that $E_i$ is big. 

This completes the proof of Lemma~\ref{small and big are alternate}
\end{proof}

Now we show that we can naturally obtain an untelescoping from 
the decomposition of $M$ by $M_{j,i}$'s. 
We explain this by giving concrete descriptions for one example. 
Giving general descriptions will easily follows from this example. 
Let $H_1 \cup_P H_2$ be a genus 4 Heegaard splitting with a 
maximal weakly reducing collection of disks 
$\Delta = \Delta_1 \cup \Delta_2$ as in Figure~4.1. 

\begin{figure}[ht]
\begin{center}
\includegraphics[width=6cm, clip]{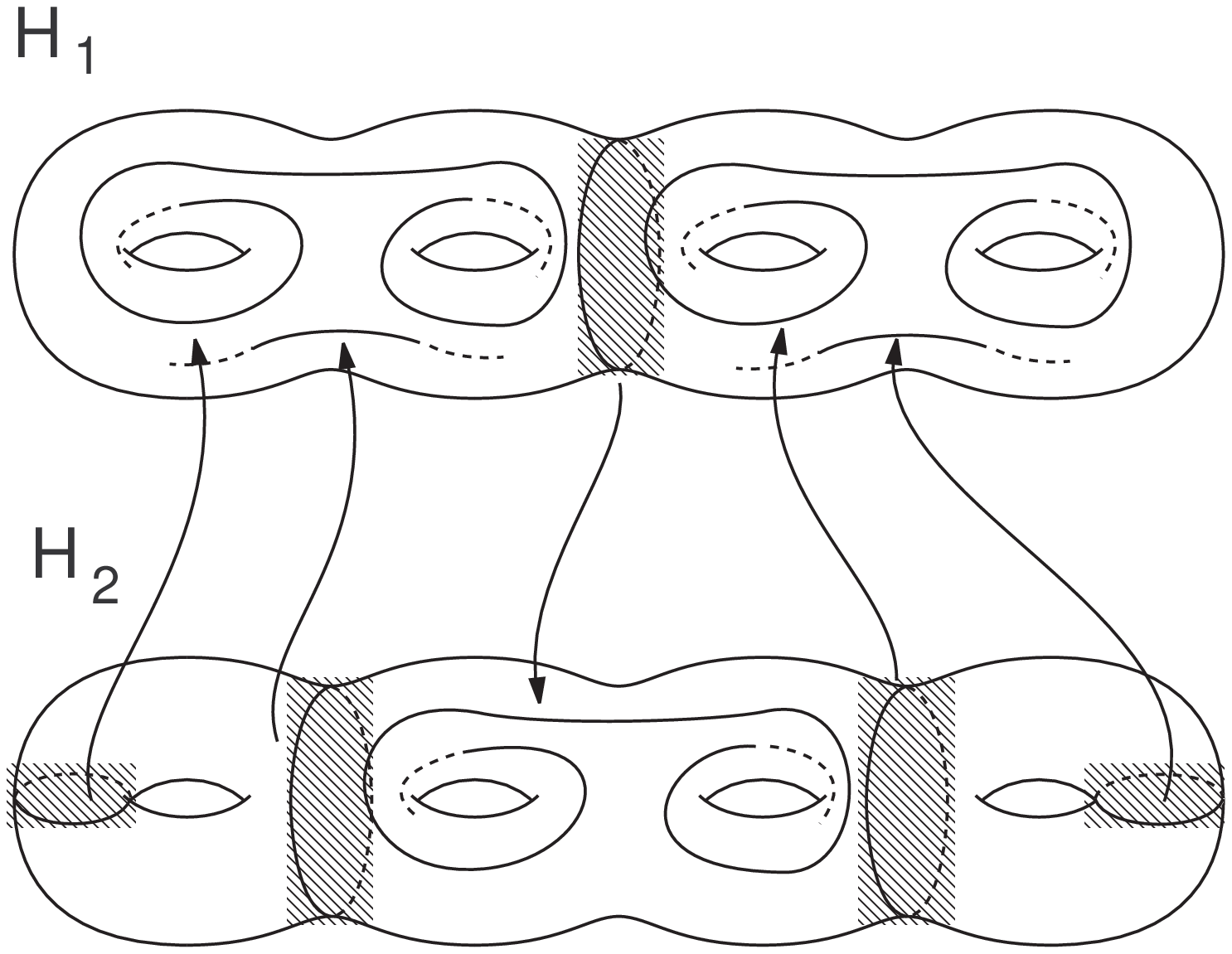}
\end{center}
\begin{center}
Figure 4.1
\end{center}
\end{figure}

In this case, 
$M$ is decomposed into three components, 
say $M_1$, $M_2$, $M_3$, 
by $\hat{P}( \Delta )$, where 
(i) $M_{1,1}$ ($M_{3,1}$ resp.) is a genus two handlebody 
which is big, and 
$M_{1,2}$ ($M_{3,2}$ resp.) is a genus 1 handlebody 
with 
$M_{1,2} \cap P$ ($M_{3,2} \cap P$ resp.) a torus with one hole, 
and 
(ii) $M_{2,2}$ is a genus two handlebody which is big, and 
$M_{2,1}$ is a 3-ball with $M_{2,1} \cap P$ an annulus. 

\begin{figure}[ht]
\begin{center}
\includegraphics[width=6cm, clip]{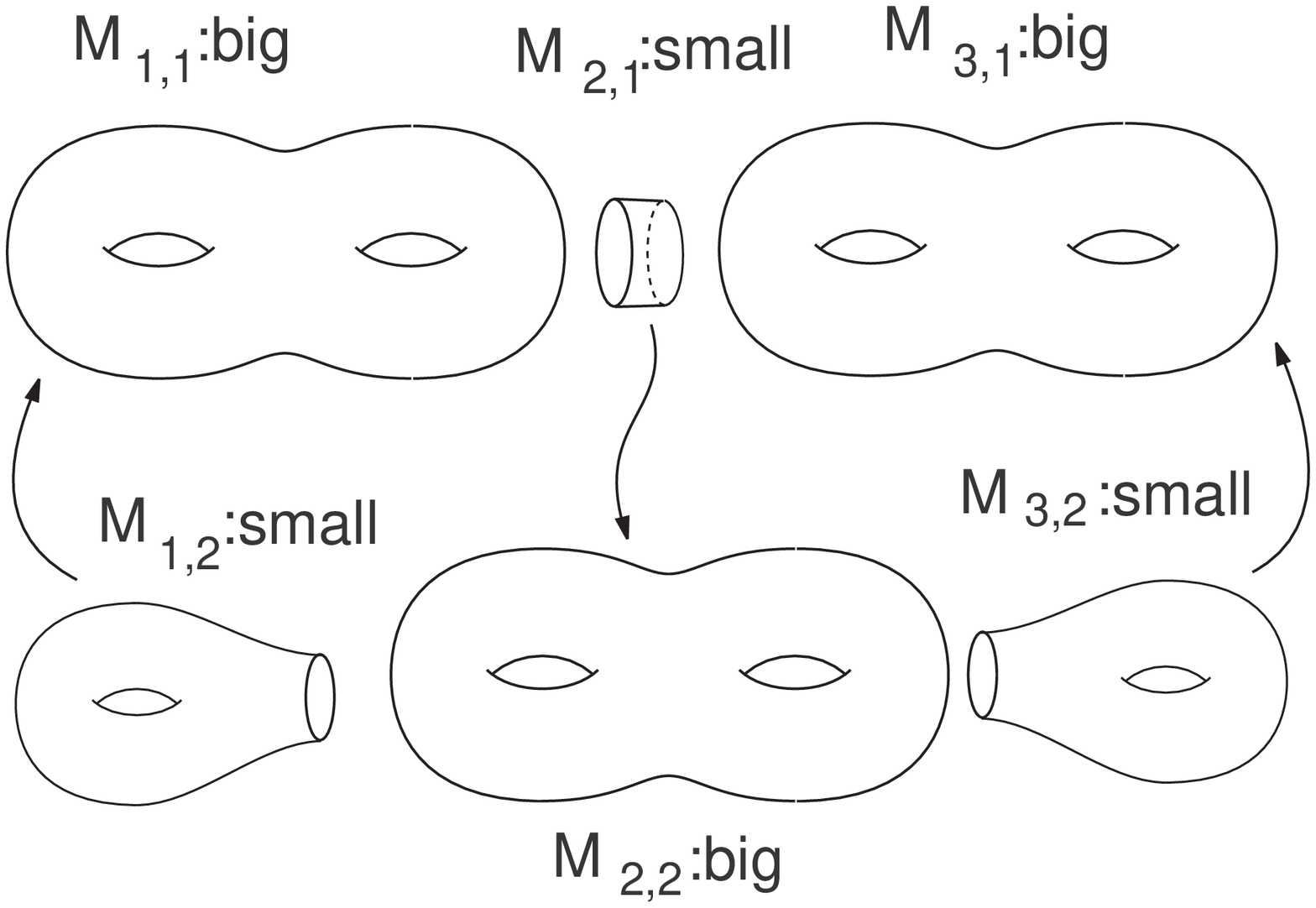}
\end{center}
\begin{center}
Figure 4.2
\end{center}
\end{figure}

Then we have the following handle decomposition of $M$. 

\medskip
$M = 
(M_{1,1} \cup M_{1,2}) \cup 
(M_{3,1} \cup M_{3,2}) \cup 
(M_{2,1} \cup M_{2,2})$

$=
\{ (0\text{-handle}) \cup 2 \times (1\text{-handle})\}
\cup
\{ 2 \times (2\text{-handle}) \cup (3\text{-handle})\}$

$\cup 
\{ (0\text{-handle}) \cup 2 \times (1\text{-handle})\}
\cup
\{ 2 \times (2\text{-handle}) \cup (3\text{-handle})\}$

$\cup 
\{ 1\text{-handle} \}
\cup
\{ 2 \times (2\text{-handle}) \cup (3\text{-handle})\}$

\medskip
Hence we can obtain an untelescoping 

$$M = 
( C_1^{(1)} \cup C_2^{(1)}) \cup 
( C_1^{(2)} \cup C_2^{(2)}) \cup 
( C_1^{(3)} \cup C_2^{(3)})$$ 

\noindent
corresponding to this handle decomposition. 

\begin{definition}\label{C-G untelescoping}
Let $\Delta$ be a weakly reducing collection of disks 
for $P$ such that each component of 
$\hat{P}( \Delta )$ is an incompressible surface in $M$, 
which is not a 2-sphere. 
We call the untelescoping of $H_1 \cup_P H_2$ obtained as above 
with such $\Delta$ a 
{\it Casson-Gordon untelescoping} 
(or {\it C-G untelescoping}) of $M$. 
\end{definition}

\begin{remark}\label{C-G induces S-T}
Let 
$( C_1^{(1)} \cup C_2^{(1)}) \cup 
\cdots 
\cup ( C_1^{(n)} \cup C_2^{(n)})$ 
be a C-G untelescoping of $H_1 \cup_P H_2$, and 
$R_k = C_1^{(k)} \cup C_2^{(k)}$ 
$(k = 1, \dots , n)$. 
By definition, $\partial R_k$ is incompressible in $R_k$. 
Hence, by Remark~\ref{Remark of S-T untelescoping}, 
we see that $C_1^{(k)} \cup C_2^{(k)}$ admits a S-T untelescoping. 
This shows that a S-T untelescoping can be regarded as a 
(possibly trivial) refinement of C-G untelescoping. 
\end{remark}

\section{Heegaard genus two link not admitting unknotting tunnel}

The {\it tunnel number} $t(L)$ of a link $L$ in the 3-sphere 
$S^3$ is the minimal number of the components of the union of mutually 
disjoint arcs, called {\it tunnels}, $\tau$ such that 
$\partial \tau \subset L$, 
and 
$\text{cl}(S^3 - N( L \cup \tau ))$ 
is a handlebody. 
We note that the definiton implies that the exteior 
$E(L)$ 
$= \text{cl}(S^3 - N(L))$ admits a Heegaard splitting of genus 
$(t(L)+1)$. 
Hence we see that the Heegaard genus of $E(L)$ is less than or equal 
to $t(L)+1$, where Heegaard genus of a 3-manifold $M$ is the 
minimal genus of the Heegaard splittings of $M$. 
Note that if we restrict our attention to Heegaard splittings 
of $(E(L), \partial E(L), \emptyset )$, then the Heegaard genus of $E(L)$ 
is exactly $t(L)+1$. 
However, if we change the partition of $\partial E(L)$, 
then they may be different. 
In this section, we give a concrete example of 
a link not satisfying the equality. 

In the remainder of this section, 
let $L = L_1 \cup L_2$ be a link as in Figure~5.1, 
where $L_1$ is a $(4,3)$ torus knot, 
and $L_2$ is a push out of a meridian curve of $L_1$, 
i.e., 
$L$ is a connected sum of $(4,3)$ torus knot and a Hopf link. 

\begin{figure}[ht]
\begin{center}
\includegraphics[width=4cm, clip]{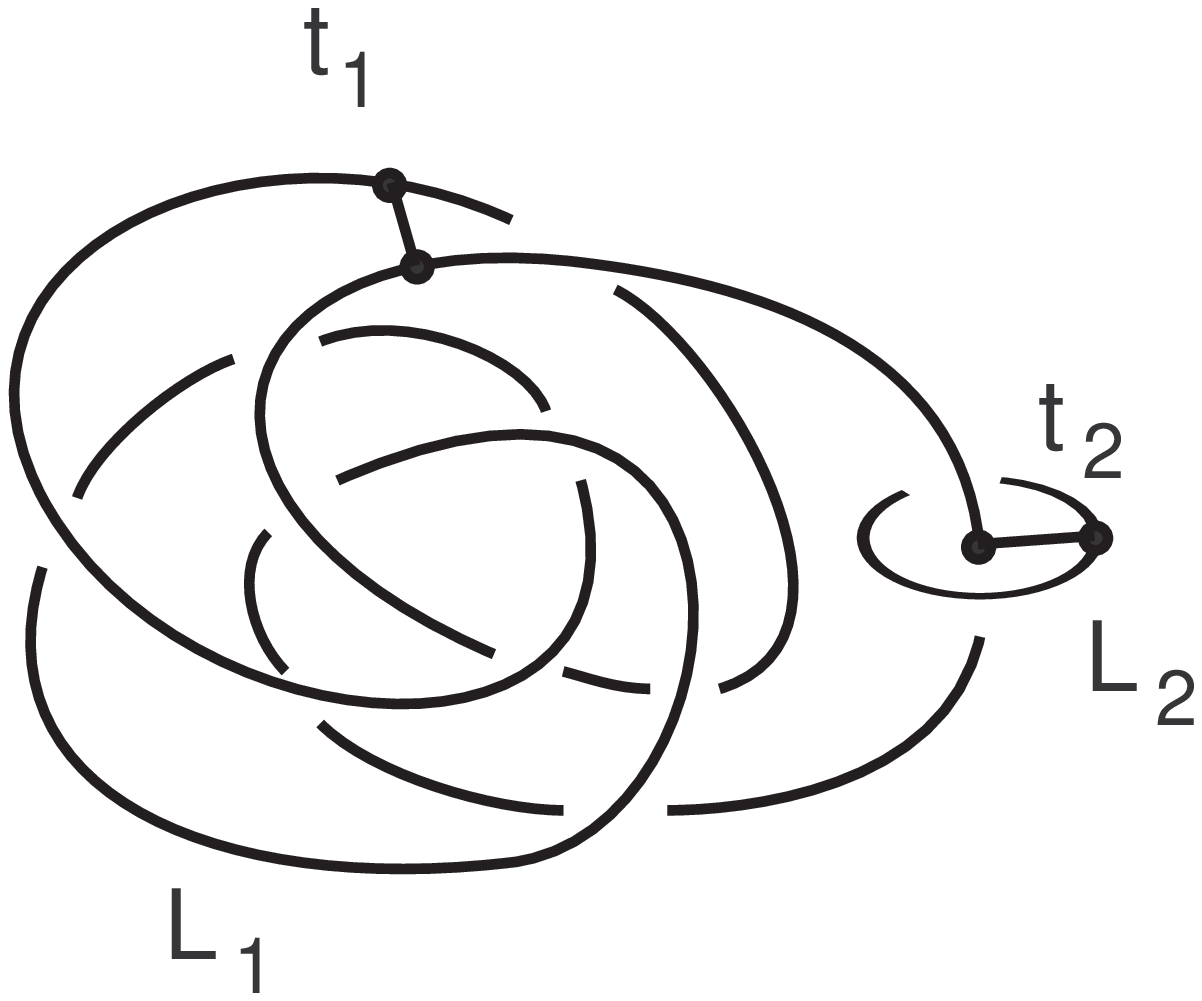}
\end{center}
\begin{center}
Figure 5.1
\end{center}
\end{figure}

\begin{proposition}\label{t(L)=2}
$t(L)=2$. 
\end{proposition}

\begin{proof}
Let $t_1$, $t_2$ be arcs as in Figure~5.1.

It is easily verified that 
$\text{cl}(S^3 - N(L \cup t_1 \cup t_2))$ 
is a genus three handlebody. 
Hence, we have $t(L) \le 2$. 

By Morimoto \cite{Mo}, it is shown that the set of two component 
composite tunnel number one links coinsides with the set of 
links each element of which is a connected sum of 
a two bridge knot and a Hopf link. 
Since $(4,3)$ torus knot is a 3-bridge knot, we see that 
$t(L) > 1$. 

Hence we have $t(L) =2$. 
\end{proof}

\begin{proposition}\label{H genus of E(L) is 2}
The Heegaard genus of $E(L)$ is 2.
\end{proposition}

\begin{proof}
It is easy to see that $E(L)$ does not admit a genus one Heegaard splitting. 
Hence it is enough to show that $E(L)$ admits a genus 2 Heegaard splitting. 

\begin{figure}[ht]
\begin{center}
\includegraphics[width=6cm, clip]{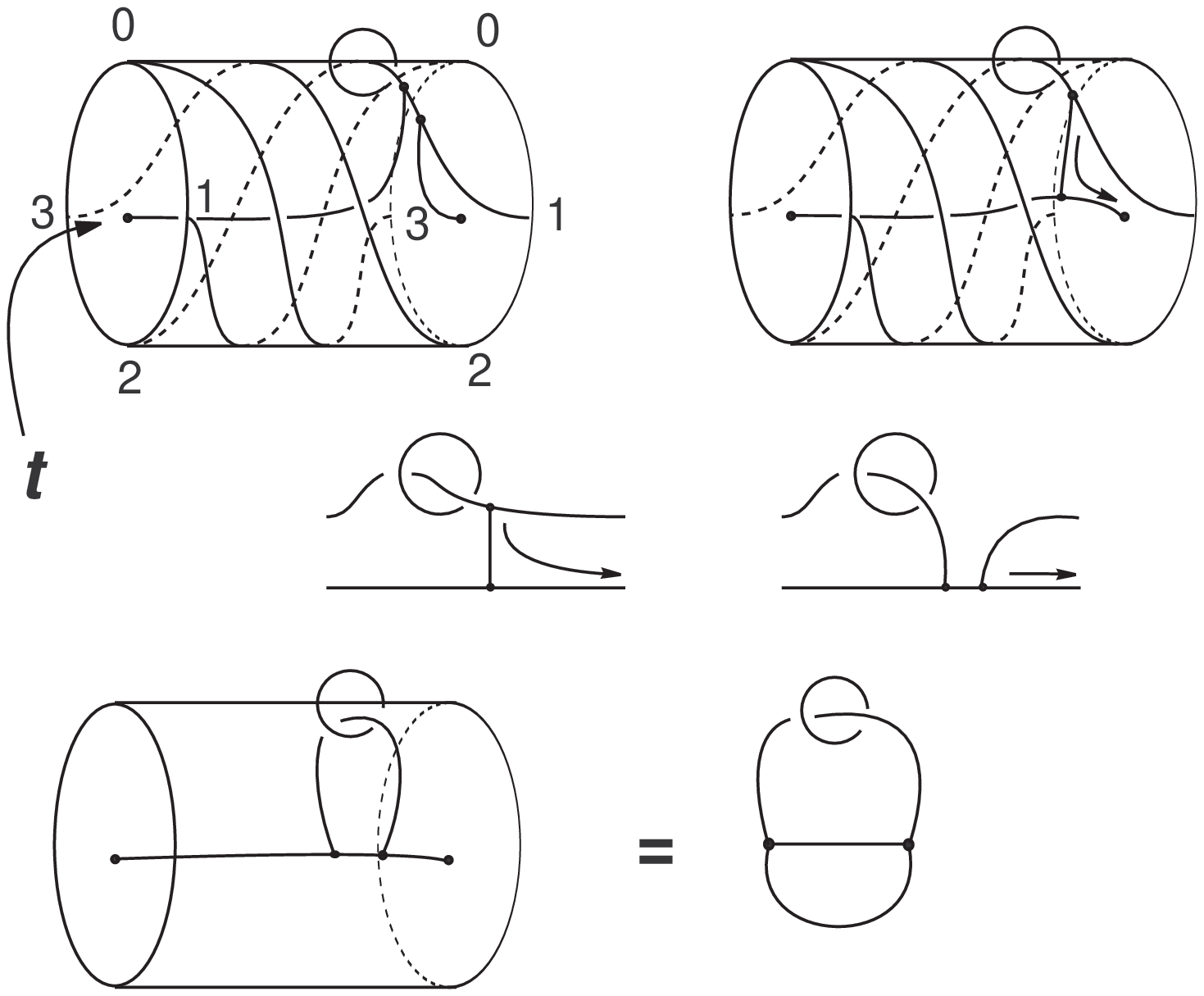}
\end{center}
\begin{center}
Figure 5.2
\end{center}
\end{figure}

Let $t$ be an arc as in Figure~5.2. 
Let $E(L_2) = \text{cl}(S^3 - N(L_2))$. 
We may suppose that $N(L_1 \cup t ) \subset \text{Int}E(L_2)$. 
The deformations in Figure~5.2 show that 
$\text{cl} (E(L_2) - N( L_1 \cup t))$ 
is a genus 2 compression body, say $C_2$, 
with $\partial_- C_2 = \partial N(L_2)$. 
Let 
$N_1 = N(L_1, N(L_1 \cup t))$, and 
$C_1 = \text{cl} (N(L_1 \cup t)-N_1)$. 
We note that $C_1$ is a genus 2 compression body with 
$\partial_- C_1 = \partial N_1$. 
Hence $C_1 \cup C_2$ gives a 
genus 2 Heegaard splitting of 
$\text{cl}(S^3 - (N_1 \cup N(L_2)))$, 
a exterior of $L$. 
\end{proof}

\section{Proof of Theorem}

Let $D(2)$ be a Seifert fibered manifold with orbit manifold a disk with 
two exceptional fibers. 
Let $L_W$ be a Whitehead link (Figure~6.1), 
$W = \text{cl}(S^3 - N(L_W))$, and 
$T_1$, $T_2$ the boundary components of $W$. 

\begin{figure}[ht]
\begin{center}
\includegraphics[width=4cm, clip]{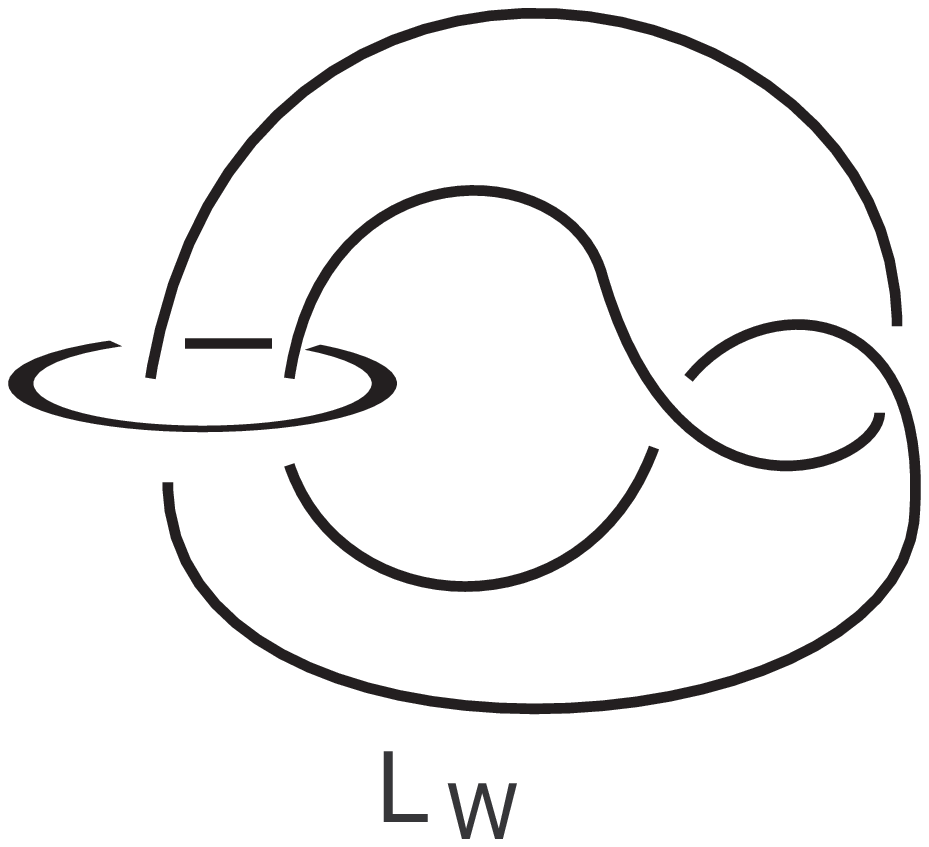}
\end{center}
\begin{center}
Figure 6.1
\end{center}
\end{figure}

Let $A(1)$ be a Seifert fibered manifold with orbit manifold 
an annulus with one exceptional fiber. 
Then let $N$ be a 3-manifold obtained from $D(2)$, $W$, and $A(1)$ 
by identifying $\partial D(2)$ and $T_1$ by a homeomorphinsm 
taking a regular fiber of $D(2)$ to a meridian curve, 
then identifying a component of $\partial A(1)$ and $T_2$ by 
a homeomorphism taking a regular fiber to a meridian curve. 
We note that $N$ is what is called a full Haken manifold 
in \cite{Ko1}, and this implies: 

\begin{proposition}\label{HG of N is 2}
The Heegaard genus of $N$ is 2. 
\end{proposition}

By Proposition~\ref{HG of N is 2}, 
we see that $N$ admits a handle decomposition 
as in the following. 

\begin{quote}
(N1) $(0\text{-handle}) \cup 2 \times (1\text{-handle}) \cup (2\text{-handle})$, 
\end{quote}

or, dually, 

\begin{quote}
(N2) $(\text{torus}\times [0,1] ) 
\cup (1\text{-handle}) 
\cup 2 \times (2\text{-handle})
\cup (3\text{-handle})$. 
\end{quote}

Let 
$N^i$ $(i=1,2)$ 
be a copy of $N$, and 
$D(2)^i$, $W^i$, $A(1)^i$, $T_1^i$, $T_2^i$ 
the pieces of $N^i$ corresponding to 
$D(2)$, $W$, $A(1)$, $T_1$, $T_2$ 
respectively. 
Let $E=\text{cl}(S^3-N(L))$, where $L = L_1 \cup L_2$ is the link in section~5. 
Note that, by the proof of Proposition~\ref{H genus of E(L) is 2}, 
$E$ admits a handle decomposition as in the following. 

\begin{quote}
(E1) 
$(\partial N(L_1) \times [0,1]) \cup (1\text{-handle}) \cup (2\text{-handle})$, 
\end{quote}

or, dually, 

\begin{quote}
(E2) 
$(\partial N(L_2) \times [0,1]) \cup (1\text{-handle}) \cup (2\text{-handle})$. 
\end{quote}

Here we note that $E$ admits a decomposition 
$E = E(4,3) \cup R$, 
where 
$E(4,3)$ is a exterior of $(4,3)$ torus knot, and 
$R$ a Seifet fibered manifold with orbit manifold a disk with two holes 
and no exceptional fibers 
(i.e., 
$R = \text{(disk with two holes)} \times S^1$), 
where a regular fiber of $R$ is identified with a meridian curve. 
Let $M$ be a closed 3-manifold obtained from 
$N_1 \cup N_2$, and $E$ by identifying their boundaries so that the 
Seifert fibrations in $A(1)^1$, $A(1)^2$, and $R$ 
do not meet on each glueing torus. 

Let $X$ be a closed Haken manifold. 
Then, by \cite{Ja}, there is a maximal perfectly embedded 
Seifert fibered manifold $\Sigma$ which is called 
a characteristic Seifert pair for $X$. 
Note that $\partial \Sigma$ consists of tori in $X$. 
If there is a pair of components of $\partial \Sigma$ 
which are parallel in $X$, then we eliminate one of them from 
the system. 
By repeating this procedure, we finally obtain a system of 
tori, say ${\cal T}$, in $X$, the elements of which are mutually 
non-parallel in $X$. 
In this paper, we call the decomposition of $X$ by ${\cal T}$ 
a {\it torus decomposition} of $X$. 
Then let $G_{\cal T}$ be the graph such that the vertices of $G_{\cal T}$ 
correspond to the components of $X-{\cal T}$, and 
the edges of $G_{\cal T}$ correspond to the components of ${\cal T}$. 
We call $G_{\cal T}$ a {\it characteristic graph} of $X$.

By the construction we immediately see that the decomposition 

$$D(2)^1 \cup_{T_1^1} W^1 \cup_{T_2^1} A(1)^1 \cup 
(R \cup E(4,3)) \cup 
A(1)^2 \cup_{T_2^2} W^2 \cup_{T_1^2} D(2)^2$$

\noindent
is a torus decomposition of $M$, where the characteristic graph is 
as follows. 

\begin{figure}[ht]
\begin{center}
\includegraphics[width=6cm, clip]{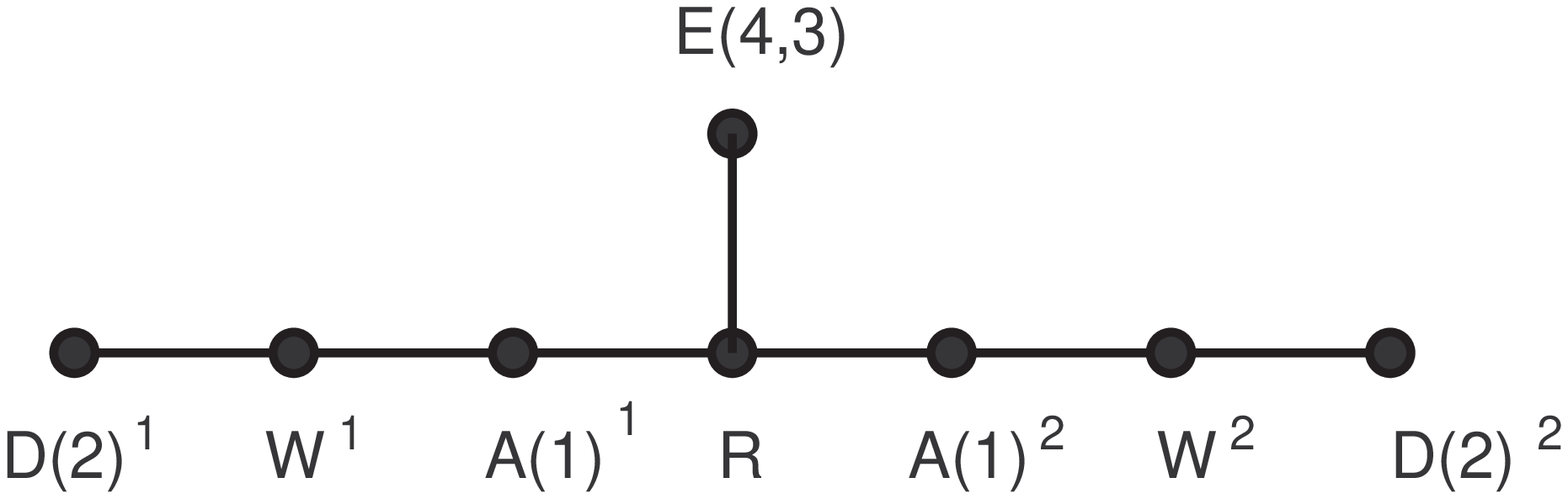}
\end{center}
\end{figure}

\begin{proposition}\label{HG of M is 4}
The Heegaard genus of $M$ is 4. 
\end{proposition}

\begin{proof}
Recall the decomposition 
$M=N_1 \cup_{\partial N(L_1)} E \cup_{\partial N(L_2)} N_2$. 
By taking the handle decompositions 
(N1) for $N_1$, 
(E1) for $E$, and 
(N2) for $N_2$, 
we see that $M$ admits a handle decoposition with 
one 0-handle, four 1-handles, four 2-handles, and one 3-handle. 
This shows that $M$ admits a genus 4 Heegaard splitting. 
Hence $g(M) \le 4$. 

It is known \cite{Ko1} 
that any Haken manifolds admitting genus $g$ Heegaard splittings 
are decomposed into at most $3g-3$ components by torus decomposition. 
Hence any Haken manifold with genus 3 Heegaard splitting is decomposed 
into at most 6 pieces by the torus decomposition. 
Recall that $M$ is decomposed into 8 pieces by the torus decomposition. 
Hence $g(M) \ge 4$. 

These show that $g(M) = 4$. 
\end{proof}

Let $V \cup_P W$ be the genus 4 Heegaard splitting obtained 
in the proof of Proposition~\ref{HG of M is 4}. 
Then the handle decopositons used there 
(
(N1) for $N_1$, 
(E1) for $E$, and 
(N2) for $N_2$
)
naturally gives an untelescoping 

$$M= 
(V_1 \cup W_1) \cup
(V_2 \cup W_2) \cup
(V_3 \cup W_3),$$

\noindent
where 
$V_1 \cup W_1$ is a genus 2 Heegaard splitting of $N_1$, 
$V_2 \cup W_2$ is a genus 2 Heegaard splitting of $E$, and 
$V_3 \cup W_3$ is a genus 2 Heegaard splitting of $N_2$. 

\begin{proposition}\label{S-T untelescoping of M}
The untelescoping 
$M= 
(V_1 \cup W_1) \cup
(V_2 \cup W_2) \cup
(V_3 \cup W_3)$ 
is a S-T untelescoping. 
\end{proposition}

\begin{proof}
In general, 
by the arguments in section~4, we easily see 
that if a genus two Heegaard splitting is 
weakly reducible, then the ambient 3-manifold is either reducible or 
admits a genus 1 Heegaard splitting. 
Each of the manifolds $N_1$, $E$, $N_2$ is not reducible or does not 
admit genus one Heegaard splitting. 
Hence we see that each $V_i \cup W_i$ $(i=1,2,3)$ is strongly irreducible. 
\end{proof}

\begin{proposition}\label{CG of M is 2 HS}
$M$ cannot be decomposed into more than two pieces by any C-G untelescoping on 
$V \cup W$. 
\end{proposition}

\begin{proof}
Assume that $M$ is decomposed into three pieces 
by a C-G untelescoping on $V \cup W$. 
Let $\Delta = \Delta_V \cup \Delta_W$ be the system of weakly reducing 
pair of disks, and 
$$M= 
(V_1 \cup W_1) \cup
(V_2 \cup W_2) \cup
(V_3 \cup W_3)$$
the C-G untelescoping. 

Let $M_i = V_i \cup W_i$ $(i=1,2,3)$. 
Recall, from section~4, that each $M_i$ has exactly one big component. 
Without loss of generality, we may suppose that $V$ contains two big components. 
Since $M$ is irreducible, 
we see that  each big comonent is a handlebody whose genus is at least two. 
It is easy to see that these together with 
Lemma~\ref{small and big are alternate} imply: 

\begin{quote}
$\Delta_V$ consists of a disk which separates $V$ into two 
genus two handlebodies which are big components. 
\end{quote}

By exchanging subscripts, if necessary, we may suppose that 
these big components correspond to $V_1$ and $V_3$. 
This implies:

\medskip
\noindent
Claim 1. Each component of $\partial M_1$, $\partial M_3$ 
is a torus, hence each component of $\partial M_2$ is a torus. 

We note that the small components of $M_2$ is just a regular 
neighborhood of $\Delta_V$ in $V$. 
Hence Claim~1 shows that the big components of $M_2$ is a genus 
two handlebody and $\partial M_2$ consists of two tori such that 
one is the boundary of $M_1$, and the other is the boundary of $M_2$. 
From these observations, we can show that the configulation of $\Delta$ 
must be as in Figure~6.2. 

\begin{figure}[ht]
\begin{center}
\includegraphics[width=6cm, clip]{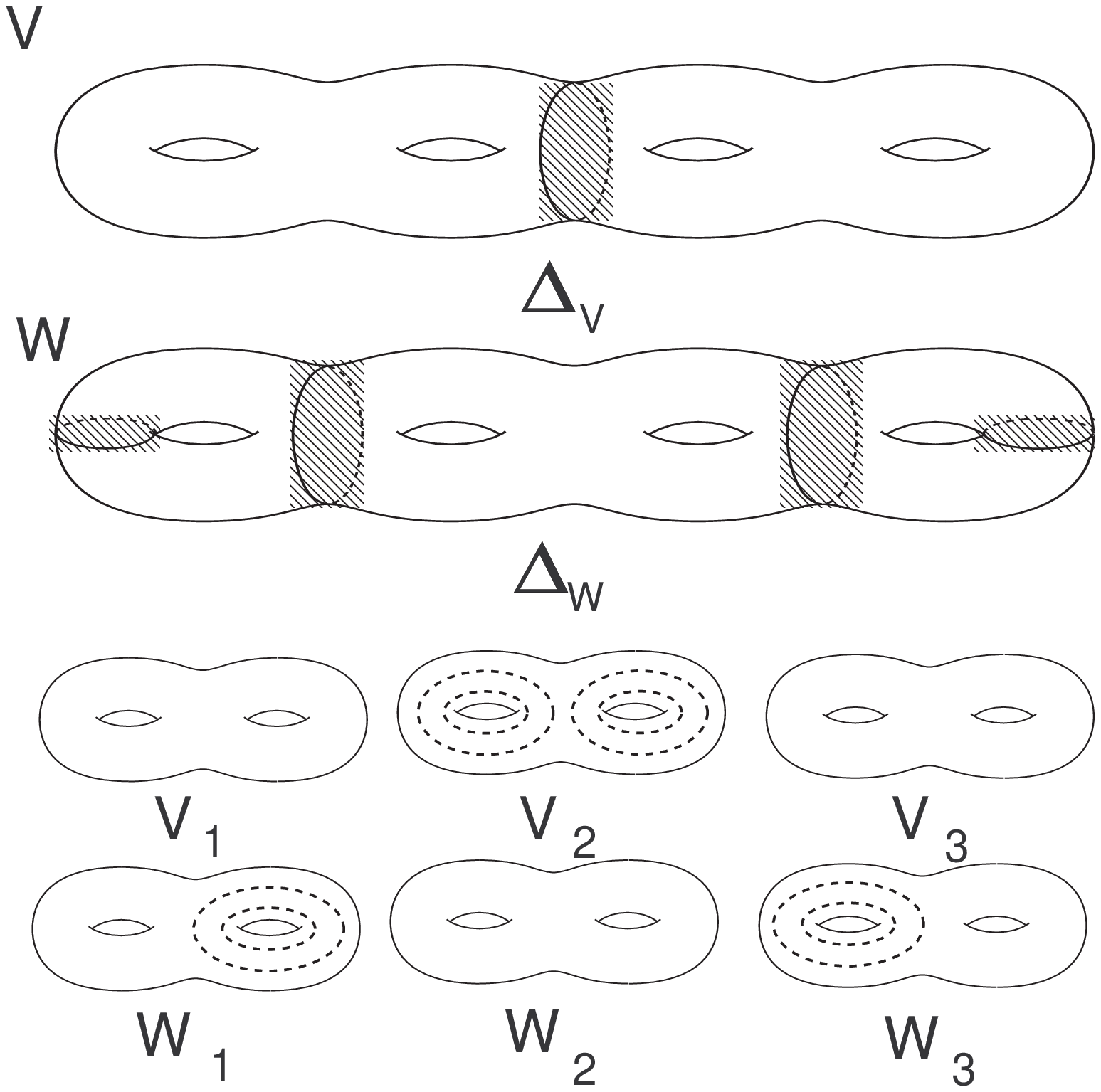}
\end{center}
\begin{center}
Figure 6.2
\end{center}
\end{figure}

\medskip
\noindent
Claim 2. 
Either $M_1 = N_1$, $M_1 = N_2$, $M_3 = N_1$, or $M_3 = N_2$. 

\begin{proof}
Note that each piece of the torus decomposition of $M$ is simple. 
Hence each incompressible torus is isotopic to a member of the 
tori giving the torus decomposition of $M$. 
Since each $M_i$ admits genus two Heegaard splitting, it 
can be decomposed into at most three pieces by the torus decomposition 
\cite{Ko1}. 
These imply that either $M_1$ or $M_3$ is decomposed into exactly three 
pieces. 
Without loss of gnerality, we may suppose that 
$M_1$ is decomposed into three pieces. 
Then note that $\partial M_1$ consists of a torus. 
These together with the examination of the charcteristic graph 
show that 
$M_1$ is either 
$D(2)^1 \cup W^1 \cup A(1)^1$, 
or 
$D(2)^2 \cup W^2 \cup A(1)^2$, 
and this proves Claim~2. 
\end{proof}

Since the argument is symmetric, we may suppose 
$M_1 = N_1$ in the remeinder of this paper. 

\medskip
\noindent
Claim 3. 
$M_3 = N_2$. 

\begin{proof}
Suppose not. 
Since $M_2$, $M_3$ are decomposed into at most three pieces 
by the torus decomposition, 
we see that 
$M_2 = E(4,3) \cup R \cup A(1)$. 
However, this contaradicts Theorem of \cite{Ko1}, 
since $R$ is not a 2-bridge link exterior. 
\end{proof}

By Claim~3, we see that $M_2 = E(L)$. 
Then, by Figure 6.2, we have $t(L) = 1$, 
contradicting Proposition~\ref{t(L)=2}. 

Finally, 
by the arguments of the proof of Proposition~\ref{S-T untelescoping of M}, 
we see that $M$ cannot be decomposed into more than three pieces, and this 
completes the proof of Proposition~\ref{CG of M is 2 HS}. 
\end{proof}

Proposition~\ref{S-T untelescoping of M} gives conclusion~1 of Theorem. 
Proposition~\ref{CG of M is 2 HS} gives conclusion~2 of Theorem. 
By the construction, we immediately have conclusion~3 of Theorem.

This completes the proof of Theorem.

\end{document}